\documentclass[11pt]{amsart}

\newtheorem{theorem}{Theorem}[section]
\newtheorem{lemma}[theorem]{Lemma}
\newtheorem{corollary}[theorem]{Corollary}
\theoremstyle{definition}

\theoremstyle{remark}

\numberwithin{equation}{section}

\newcommand{\Rls}{{\mathbb R}}

\begin{document}               

\title[pinching estimates]
{Pinching estimates and motion of hypersurfaces by curvature functions}
\author{Ben Andrews} 
\address{CMA, ANU, ACT 0200, Australia}
\email{andrews@maths.anu.edu.au}
\thanks{Research supported by a grant from the Australian Research Council}
\subjclass{Primary 53C44; Secondary  35K55, 35B50.}
\begin{abstract}
Second derivative pinching estimates are proved for a class of elliptic and 
parabolic equations, including motion of hypersurfaces by curvature functions such
as quotients of elementary symmetric functions of curvature.  The estimates imply
convergence of convex hypersurfaces to spheres under these flows, improving earlier
results of B.~Chow and the author.  The result is obtained via a detailed analysis
of gradient terms in the equations satisfied by second derivatives.
\end{abstract}
\maketitle

\section{Introduction}

The aim of this paper is to provide some insights into second-derivative
estimates for fully nonlinear elliptic and parabolic equations.  In particular, the
paper will explore the nonlinear terms which arise in the equations 
satisfied by second derivatives of solutions, and 
introduce some tools for understanding these.  The result is
a new pinching estimate for second derivatives, which improves several
previously known results.  The estimate has a number of applications, but this paper 
emphasises particularly  the implications for the behaviour of convex 
hypersurfaces moving with speed given by a function of principal curvatures.

\subsection{Second derivative estimates}

To motivate the main result, consider as a model problem fully nonlinear 
scalar parabolic equations of the form
\begin{equation}\label{eq:FNPE}
\frac{\partial u}{\partial t} = F(D^2u)
\end{equation}
where $F$ is a smooth function of the components of $D^2u$.

Differentiation of Equation \eqref{eq:FNPE} yields a very nice system of equations 
satisfied by the first derivatives of a solution:
\begin{equation}\label{eq:sysdu}
\frac{\partial}{\partial t}D_ju = \dot F^{kl}D_kD_l\left(D_ju\right)
\end{equation}
where we sum over repeated indices, and  $\{\dot F^{kl}\}$ is the matrix of partial 
derivatives of $F$ with respect to the components of its argument:
\[
\frac{\partial}{\partial s}F(A+sB)\Big|_{s=0} = \dot F^{kl}(A)B_{kl}.
\]
However, after differentiating once more, the resulting system of equations for the
second derivatives is more complicated:
\begin{equation}\label{eq:sysd2u}
\frac{\partial}{\partial t}\left(D_iD_ju\right)
= \dot F^{kl}D_kD_l\left(D_iD_ju\right)
+ \ddot F^{kl,pq} D_iD_kD_lu D_jD_pD_qu
\end{equation}
where $\ddot F$ is the second derivative of $F$:
\[
\frac{\partial^2}{\partial s^2}F(A+sB)\Big|_{s=0} = \ddot F^{kl,pq}(A)B_{kl}B_{pq}.
\]
The second term on the right-hand side of the system \eqref{eq:sysd2u} is an 
obstacle to simple applications of the maximum principle to control the behaviour
of second derivatives of solutions, since it is difficult to obtain useful
information on its sign.   This applies both to arguments using the
classical maximum principle and to those using other tools such as the
Aleksandrov-Bakelman maximum principle.  I will restrict the discussion here to the
classical setting.  

The main result of the paper is essentially the following:  If $F$ is concave
as a function of the second derivatives, and also `inverse-concave', meaning that the 
function $F^*$ defined by $F^*(A)=-F(A^{-1})$ is concave, then the ratio of minimum
eigenvalue to trace of $D^2u$ never decreases below its initial minimum.
The precise statement is given in Theorem \ref{thm:main}.     
The proof requires a detailed understanding of the nonlinear terms arising in Equation
\eqref{eq:sysd2u}, and includes several useful tools for understanding these.  Also 
important in the application of the main result is a new maximum principle for
tensors which is given in Theorem \ref{thm:newMP}.

\subsection{Evolving hypersurfaces}

There has been considerable previous work on convex hypersurfaces moving by 
curvature flows, and the most relevant here is where the speed function is a
homogeneous degree one, monotone increasing function of the principal 
curvatures.  The first such flow considered was the flow by mean curvature, which
was treated by G.~Huisken \cite{Hu1}.  He proved that convex hypersurfaces contract
to points in finite time under this flow, with spherical limiting shape.  B.~Chow
proved a similar result  for the motion of an $n$-dimensional hypersurface by the
$n$th root of the Gauss curvature \cite{Ch1}.   He also proved a result for
motion by the square root of the scalar curvature \cite{Ch2}, but in that case
a stronger assumption than convexity was required for the initial hypersurface. 
The author considered a very general class of homogeneous degree one flows in
\cite{A1}, and proved the general result if the  speed is a convex function of
principal curvatures (as is the mean curvature) or if it is concave in the
principal curvatures and vanishes when any principal curvature approaches
zero (as in the case of the $n$th root of the Gauss curvature).  More
generally, if the speed is concave in the principal curvatures, the result
holds as long as we assume a strong enough pinching condition on the initial
hypersurface.  The pinching estimate proved in this paper is aimed at
removing the latter restriction for a wide class of flows of interest,
including the flow by square root of scalar curvature treated by Chow.  

The new result for contracting hypersurfaces is as follows:   

\begin{theorem}\label{thm:contraction}
Let $f$ be a smooth symmetric  function defined on the positive cone 
$\Gamma_+=\{(x_1,\dots,x_n)\in{\mathbb R}^n:\ x_i>0,\ i=1,\dots,n\}$ in ${\mathbb
R}^n$, $n\geq 2$, which is homogeneous of degree one and strictly monotone 
increasing in each
argument.   Suppose that either
\begin{enumerate}
\item $n=2$, or
\item $f$ is convex, or 
\item $f$ is concave on $\Gamma_+$ and zero on the boundary of $\Gamma_+$, or 
\item both $f$ and the function $f^*$ given by 
$f^*(x_1,\dots,x_n)=-f(x_1^{-1},\dots,x_n^{-1})$ are concave on $\Gamma_+$.
\end{enumerate}

Let $x_0: M^n\to{\mathbb R}^{n+1}$ be a smooth, strictly convex embedding.  Then 
there exists a unique maximally extended solution $x: M\times [0,T)\to{\mathbb
R}^{n+1}$ of 
\begin{align}\label{eq:flow}
\frac{\partial x(p,t)}{\partial t}& = -f(\kappa_1(p,t),\dots,\kappa_n(p,t))\nu(p,t),
\quad (p,t)\in M\times[0,T);\\
x(p,0)&=x_0(p),\quad p\in M,\notag
\end{align}
where $\kappa_1(p,t),\dots,\kappa_n(p,t)$ are the principal curvatures of the 
embedding $x_t(.)=x(.,t)$ at the point $x(p,t)$, and $\nu(p,t)$ is the
outward-pointing unit normal vector to $x_t(M)$ at $x(p,t)$.  The map $x_t$ converges 
to a constant $z\in{\mathbb R}^{n+1}$ as $t$ approaches $T$, and the rescaled
embeddings
$\tilde x_t = \frac{x_t-z}{\sqrt{T-t}}$ converge in $C^{\infty}$ to a limit with 
image a sphere of radius $\sqrt{2f(1,\dots,1)}$ centred at the origin.

The same result holds if $f$ is merely concave, provided the initial embedding is 
such that 
$$
\sup_{p\in M} \frac{\kappa_1(p,0)+\dots+\kappa_n(p,0)}{f(\kappa_1(p,0),
\dots,\kappa_n(p,0))} < \liminf_{(x_1,\dots,x_n)\to\partial\Gamma_+}
\frac{x_1+\dots+x_n}{f(x_1,\dots,x_n)}.
$$
\end{theorem}

The new ingredient is the last of the four conditions allowed for $f$ (the 
first is treated in a recent paper by the author \cite{A3} making use of
some new regularity results special to two dimensions proved in \cite{A2},
and the second and third cases were proved in \cite{A1}).  

\section{A class of symmetric functions} 

The statement of Theorem \ref{thm:contraction} brings interest to a certain class
of symmetric functions defined on the positive cone $\Gamma_+$ in ${\mathbb R}^n$.
In this section I will discuss this class in some detail.

The class of interest, denoted ${\mathcal C}_n$, consists of 
functions on $\Gamma_+$ which are
\begin{itemize}
\item smooth ($C^\infty$);
\item homogeneous of degree one:  $f(cx)=cf(x)$ for $c>0$;
\item strictly monotone increasing: $\frac{\partial f}{\partial x_i}>0$ for each $i$;
\item  concave; and 
\item inverse-concave:  $f^*(x_1,\dots,x_n) = -f(x_1^{-1},\dots,x_n^{-1})$ defines a 
concave function on $\Gamma_+$.  
\end{itemize}
Also important is the subclass $\mathcal S_n$ consisting of
functions in $\mathcal C_n$ which are symmetric: 
$f(x_{\sigma(1)},\dots,x_{\sigma(n)}) = f(x_1,\dots,x_n)$ for any permutation 
$\sigma$.
These are precisely the functions which satisfy the fourth condition 
in Theorem \ref{thm:contraction}.  Note that the main result, Theorem \ref{thm:main},
does not require any homogeneity condition, and so applies to a somewhat larger
class than ${\mathcal S}_n$.

Before giving examples of functions in ${\mathcal S}_n$, I will give some 
useful methods of constructing new examples from old.

\begin{theorem}\label{thm:Snclass} A homogeneous degree one function
$f:\Gamma_+\to{\mathbb R}$ is in ${\mathcal C}_n$ if and only if the following
conditions hold everywhere on $\Gamma_+$:
\begin{enumerate}
\item $\dot f^i = 
\frac{\partial f}{\partial x_i}>0$ for each $i$;
\item $\ddot f^{ij}=\frac{\partial^2f}{\partial x_i\partial x_j}$ is a non-positive
matrix;
\item $\ddot f^{ij}+2\frac{\dot f^i}{x_i}\delta_{ij}$ is a non-negative matrix.
\end{enumerate}
\end{theorem}

\begin{proof}
The only non-trivial point is that the third condition is equivalent to the
concavity of $f^*$.  To see this, compute the derivatives of $f^*$ at
$(z_1,\dots,z_n)$, where $z_i=x_i^{-1}$:
\begin{align*}
\frac{\partial f^*}{\partial z_i} &=
\frac{\partial f}{\partial x_i}x_i^2;\\
\frac{\partial^2 f^*}{\partial z^i\partial z^j}
&=-\frac{\partial^2 f}{\partial x_i\partial x_j}x_i^2x_j^2
-2\frac{\partial f}{\partial x_i}x_i^3\delta_{ij}.
\end{align*}
Multiplying the last identity by $x_i^{-2}x_j^{-2}$ gives the result.
\end{proof}

\begin{corollary}\label{cor:mean}
The algebraic mean $H=\frac1n\sum_ix_i$ is in ${\mathcal S}_n$.
\end{corollary}

\begin{proof}
In this case $\ddot f=0$ and $\dot f>0$.
\end{proof}

\begin{theorem}\label{thm:power}
If $f\in {\mathcal C}_n$ and $r\in[-1,1]\backslash\{0\}$, then the function $f_r$
given by
$$
f_r(x_1,\dots,x_n) = \left(f(x_1^r,\dots,x_n^r)\right)^{\frac1r}
$$
is in ${\mathcal C}_n$.
\end{theorem}

\begin{proof}
Compute the first and second derivatives of $f_r$ at $z_i=x_i^{1/r}$:
\begin{align*}
\frac{\partial f_r}{\partial z_i}
&=f^{\frac1r-1}\dot f^iz_i^{r-1}>0;\\
\frac{\partial^2 f_r}{\partial z_i\partial z_j}
&=rf^{\frac1r-1}\ddot f^{ij}z_i^{r-1}z_j^{r-1}+(1-r)f^{\frac1r-2}\dot f^i\dot f^j
z_i^{r-1}z_j^{r-1}\\
&\quad\null+(r-1)f^{\frac1r-1}\dot f^iz_i^{r-2}\delta_{ij}\\
&=rf^{\frac1r-1}z_i^{r-1}z_j^{r-1}
\left(\ddot f^{ij}-\left(\frac{r-1}r\right)\frac{\dot f^i\dot f^j}{f}
+\left(\frac{r-1}r\right)\frac{\dot f^i}{x_i}\delta_{ij}\right)
\end{align*}
Since $f_r$ is homogeneous of degree one, the Euler relation implies that 
$\sum_i\frac{\partial^2f_r}{\partial z_i\partial
z_j}z_i=0$, so the radial vector is a null eigenvector.  Therefore to prove 
concavity of $f_r$ it suffices to
consider the restriction of the second derivatives to the transversal subspace
$S=\{\xi:\ \dot f^i\xi_i=0\}$.  There the bracket becomes
\[
\ddot f^{ij}+\frac{r-1}{r}\frac{\dot f^i}{x_i}\delta_{ij}.
\]
There are two cases to consider:  If $0<r\leq 1$, then $r-1\leq 0$ and the
bracket is non-positive.
If $-1\leq r<0$, then $\frac{r-1}r\geq 2$, and the bracket is non-negative by the
third point in Theorem \ref{thm:Snclass}.  In the latter case the coefficient is
negative, so in both cases $f_r$ is concave.

To establish the inequality in the third part of Theorem
\ref{thm:Snclass} for the function $f_r$, it suffices to show the stronger
inequality 
\begin{equation}\label{eq:stronginvconc}
\ddot f_r^{ij}+2\frac{\dot f_r^i}{z_i}\delta_{ij}
-2\frac{\dot f_r^i\dot f_r^j}{f_r}\geq 0.
\end{equation}
The expressions above give
\begin{align*}
\ddot f_r^{ij}&+2\frac{\dot f_r^i}{z_i}\delta_{ij}-2\frac{\dot f_r^i\dot f_r^j}{f_r}\\
&=rf^{\frac1r-1}z_i^{r-1}z_j^{r-1}
\left(\ddot f^{ij}-\frac{r+1}{r}\frac{\dot f^i\dot f^j}{f}+
\frac{r+1}{r}\frac{\dot f^i}{x_i}\delta_{ij}\right).
\end{align*}
As before the bracket has the radial vector as a null eigenvector, and on the
subspace $S$ it becomes $\ddot f^{ij}+\frac{r+1}{r}\frac{\dot f^i}{x_i}\delta_{ij}$,
which is bounded below by $\ddot f^{ij}+2\frac{\dot f^i}{x_i}\delta_{ij}$ if $0<r\leq
1$, and bounded above by $\ddot f^{ij}$ if $-1<r<0$.
\end{proof}

\begin{corollary}\label{cor:altdef}
A homogeneous degree one monotone increasing
function
$f$ on
$\Gamma_+$ is in
${\mathcal C}_n$ if and only if it is concave and the function $f_{-1}$ defined
in Theorem \ref{thm:power} is concave.
\end{corollary}

\begin{proof} Concavity of $f_{-1}$ is equivalent to the inequality \eqref{eq:stronginvconc} (with $r=1$), hence stronger than concavity of $f^*$.
\end{proof}

\begin{corollary}\label{cor:power}
The power-means $H_r=\left(\frac1n\sum_ix_i^r\right)^{\frac1r}$ are in
${\mathcal S}_n$ for $|r|\leq 1$.
\end{corollary}

\begin{proof}
By Corollary \ref{cor:mean} we can take $f=H$ in Theorem \ref{thm:power} to obtain the
result for $r\neq 0$.  The result holds also for $r=0$ since $H_r$ converges locally
uniformly to $H_0=(x_1\dots x_n)^{1/n}$ as $r\to 0$.
\end{proof}

A commonly defined class involves the elementary symmetric functions of
principal curvatures:
$$
S_k(x_1,\dots,x_n) = \frac1{{n\choose k}}\sum_{1\leq
i_1<\dots<i_k\leq n} x_{i_1}\dots x_{i_k},\quad k=1,\dots,n.
$$

\begin{theorem}
For $0\leq k<n$, $\frac{S_{k+1}}{S_k}\in{\mathcal S}_n$.
\end{theorem}

\begin{proof}
The concavity is proved in \cite[Theorem 15.16]{Lb}, and the inverse-concavity
follows since $f_{-1}=\frac{S_{n-k}}{S_{n-k-1}}$ is of the same kind.  
\end{proof}

To complete
our discussion and give a satisfyingly large class of examples, we note the following:

\begin{theorem}
If $f_1,\dots,f_k$ are in ${\mathcal C}_n$, and $\varphi\in\mathcal C_k$, then 
the function $f$ defined by
\[
f(x_1,\dots,x_n) = \varphi(f_1(x_1,\dots,x_n),\dots,f_k(x_1,\dots,x_n))
\]
is in ${\mathcal C}_n$.  If the $f_i$ are in $\mathcal S_n$, then so is $f$.
\end{theorem}

\begin{proof}
It suffices to show that $f$ is concave under these conditions, because 
\[
f_{-1}(x_1,\dots,x_n)=\varphi_{-1}((f_1)_{-1}(x_1,\dots,x_n),\dots,(f_k)_{-1}
(x_1,\dots,x_n)),
\]
and by Corollary \ref{cor:altdef}, $\varphi_{-1}\in\mathcal C_k$ and
$(f_i)_{-1}\in\mathcal C_n$ for each $i$, so $f_{-1}$ must also be concave.
The derivatives of $f$ are as follows:
\begin{align*}
\dot f^i&=\sum_p\dot\varphi^p\dot f_p^i;\\
\ddot f^{ij}&=\sum_{p,q}\ddot\varphi^{pq}\dot f_p^i\dot f_q^j+\sum_p\dot\varphi^p\ddot
f_p^{ij},
\end{align*}
so for any $\xi=(\xi_1,\dots,\xi_n)$, writing $\eta_p=\sum_i\dot f_p^i\xi_i$ for 
$p=1,\dots,k$,
\[
\ddot f(\xi,\xi) = \ddot\varphi(\eta,\eta)+\sum_k\dot\varphi^p\ddot f_p(\xi,\xi).
\]
This is non-positive since $\ddot\varphi\leq 0$, $\ddot f_p\leq 0$, and
$\dot\varphi^p>0$.
\end{proof}

In particular, the class ${\mathcal S}_n$ is closed under multiplication by
positive scalars, addition, and taking weighted geometric means.  Therefore the
following examples involving elementary symmetric functions are in ${\mathcal S}_n$:
\begin{itemize}
\item $S_k^{\frac1k}$, since this is the geometric mean of $\frac{S_{j+1}}{S_j}$
for $j=0,\dots,k-1$;
\item $\left(\frac{S_k}{S_l}\right)^{\frac1{k-l}}$ for $n\geq k>l\geq 0$, since
this is the geometric mean of $\frac{S_{j+1}}{S_j}$ for $j=l,\dots,k-1$;
\item $S_n^{\alpha_n}S_{n-1}^{\alpha_{n-1}-\alpha_n}\dots S_2^{\alpha_2-\alpha_3}
S_1^{\alpha_1-\alpha_2}$, if $\alpha_i\geq 0$ and $\sum_i\alpha_i=1$ (the
general form of a weighted geometric mean of $\{\frac{S_{j+1}}{S_j},\
j=0,\dots,n-1\}$);
\item positive linear combinations of any of the above examples.
\end{itemize}

All of the above examples can be used in Theorem \ref{thm:contraction}.  Note that
of these, relatively few are covered by the previously known results:  Of the 
functions
$S_k^{1/k}$, only $k=1$ and $k=n$ were known;
for $f=(S_k/S_l)^{1/(k-l)}$, $k>l$, only those with $k=n$ or $k=1$ were known;
and for the more general combinations 
$S_n^{\alpha_n}S_{n-1}^{\alpha_{n-1}-\alpha_n}\dots S_2^{\alpha_2-\alpha_3}
S_1^{\alpha_1-\alpha_2}$ with $\alpha_i\geq 0$ and $\sum\alpha_i=1$, 
only those with 
$\alpha_n>0$ were known (except $\alpha_1=1$).

Note also that the previously known results allowed speeds given by the power means
$H_r$ for $r\geq 1$ (convex case) or $r\leq 0$ (concave and zero on the boundary of
the positive cone).  The new result therefore extends this to cover all of the
remaining values of $r$.

\section{A maximum principle for tensors}

This section refines the following well-known result  from \cite{Ha1}:

\begin{theorem}[\cite{Ha1}, Theorem 9.1]
Let $M$ be a compact manifold with a (possibly time-dependent) Riemannian metric $g$ and Levi-Civita connection $\nabla$.  Let $S_{ij}$ be a smooth symmetric tensor field  satisfying
\[
\frac{\partial}{\partial t}S_{ij} = \Delta S_{ij} + u^k\nabla_kS_{ij} + N_{ij}
\]
on some time interval $[0,T]$, where $u$ is smooth,
and $N(p,t)(v,v)\geq 0$ whenever $v$ is a null eigenvector of $S(p,t)$.  
If $S_{ij}$ is positive definite everywhere at time 
$t=0$, then it remains so on $0\leq t\leq T$.
\end{theorem}

This generalises easily to the case where the manifold has boundary, and the tensor 
field satisfies a more general evolution equation of the form
\[
\frac{\partial}{\partial t}S_{ij} = a^{kl}\nabla_k\nabla_lS_{ij} +
u^k\nabla_kS_{ij}+N_{ij}
\]
where $a^{kl}$ is smooth and positive definite at each point and time.  
The result also remains true for connections other than the Levi-Civita 
connection.   The new result makes particular use of the latter observation:

\begin{theorem}\label{thm:newMP}
Let $S_{ij}$ be a smooth time-varying symmetric tensor field on a compact manifold 
$M$ (possibly with boundary), 
satisfying
\[
\frac{\partial}{\partial t}S_{ij} = a^{kl}\nabla_k\nabla_lS_{ij} +
u^k\nabla_kS_{ij}+N_{ij}
\]
where $a^{kl}$ and $u$ are smooth, $\nabla$ is a (possibly 
time-dependent) smooth symmetric connection, and $a^{kl}$ is positive definite
everywhere.  Suppose that
\[
N_{ij}v^iv^j +
\sup_{\Gamma}2a^{kl}\left(2\Gamma_k^p\nabla_lS_{ip}v^i
-\Gamma_k^p\Gamma_l^q
S_{pq}\right\}\geq 0
\]
whenever $S_{ij}\geq 0$ and $S_{ij}v^j=0$.
If $S_{ij}$ is positive definite
everywhere on $M$ at time $t=0$ and on $\partial M$ for $0\leq t\leq T$, 
then it is positive on $M\times[0,T]$.
\end{theorem}

The result amounts to the observation that in a parabolic equation of 
this kind for a tensor field (in contrast to the scalar case), the leading elliptic term 
can be squeezed to yield a non-trivial extra term.  In many situations this extra term is
not useful, but in equations like \eqref{eq:sysd2u} the `reaction' term $N_{ij}$ is quadratic in the derivatives of $S$, and the extra term given by Theorem \ref{thm:newMP} is also of this kind if $\Gamma$ is taken to be linear in the derivatives of $S$.  The extra term that results is crucial in the
proof of Theorem \ref{thm:main}.

\begin{proof}
The tensor inequality $S_{ij}\geq 0$ is equivalent to an inequality for
a function on the tangent bundle of $M$:
\[
Z(p,v) = S(p)(v,v)\geq 0
\]
for all $p\in M$ and $v$ in $T_pM$.  Let $p$ be a point where $S(p)$ has a null 
eigenvector $v$.  Choose coordinates $x^1,\dots,x^n$ for $M$ near $p$ such that
the connection coefficients of $\nabla$ vanish at $p$.  Then any vector in $T_qM$
for $q$ near $p$ has the form $\sum_{i=1}^n\dot x^i\partial_i$, so $TM$ is
described locally by the $2n$ coordinates $x^1,\dots,x^n$ and $\dot x^1,\dots,\dot
x^n$.  The coordinates can be chosen so that $v=\partial_1$.

At $(p,v)$ the first derivatives of $Z$ must vanish, so that
\begin{equation}\label{eq:xderiv}
0=\frac{\partial Z}{\partial x^k} = \left(\frac{\partial}{\partial x^k} S_{ij}\right)\dot x^i\dot x^j=
\frac{\partial}{\partial x^k}S_{11},
\end{equation}
and
\begin{equation}\label{eq:vderiv}
0=\frac{\partial Z}{\partial \dot x^k}=2S_{ij}\dot x^i\delta_k^j=2S_{k1},
\end{equation}
for $k=1,\dots,n$.  Equation \eqref{eq:vderiv} implies  that $v$ is a null eigenvector of $S(p)$.

Now consider the second order conditions implied by minimality:
The second derivatives form a $(2n)\times(2n)$ matrix which is non-negative
at $(p,v)$.  The entries in this matrix are as follows:
\begin{align*}
\frac{\partial^2Z}{\partial x^k\partial x^l}&=
\frac{\partial^2}{\partial x^k\partial x^l}S_{11};\\
\frac{\partial^2Z}{\partial x^k\partial\dot x^q}
&=2\frac{\partial}{\partial x^k}S_{1q};\\
\frac{\partial^2Z}{\partial\dot x^p\partial\dot x^q}
&=2S_{pq}.
\end{align*}
For any $\Gamma$ this implies the inequality
\begin{align}\label{eq:minineq}
0&\leq a^{kl}\left(\frac{\partial^2Z}{\partial x^k\partial x^l}
-\Gamma_k^p\frac{\partial^2Z}{\partial\dot x^p\partial x^l}
-\Gamma_l^q\frac{\partial^2Z}{\partial x^k\partial\dot x^q}
+\Gamma_k^p\Gamma_l^q\frac{\partial^2Z}{\partial\dot x^p\partial\dot x^q}\right)\notag\\
&=a^{kl}\left(\frac{\partial^2S_{11}}{\partial x^k\partial x^l}
-2\Gamma_k^p\frac{\partial S_{1p}}{\partial x^l}
-2\Gamma_l^q\frac{\partial S_{1q}}{\partial x^k}
+2\Gamma_k^p\Gamma_l^qS_{pq}\right)\notag\\
&=a^{kl}\left(\frac{\partial^2S_{11}}{\partial x^k\partial x^l}
-4\Gamma_k^p\frac{\partial S_{1p}}{\partial x^l}
+2\Gamma_k^p\Gamma_l^qS_{pq}\right).
\end{align}
In the coordinates chosen above, the coordinate derivatives at $p$ relate to the
covariant derivatives as follows:
\begin{equation}\label{eq:identderiv}
\nabla_k\nabla_lS_{11} = \frac{\partial^2S_{11}}{\partial x^l\partial x^l}
-2S(\nabla_k\nabla_l\partial_1,\partial_1)=\frac{\partial^2S_{11}}
{\partial x^k\partial x^l},
\end{equation}
in view of the first order condition \eqref{eq:vderiv}.

The function $Z$ satisfies the scalar evolution equation
\[
\frac{\partial Z}{\partial t} = a^{kl}\nabla_k\nabla_lS_{ij}\dot x^i\dot x^j +
u^k\nabla_kS_{ij}\dot x^i\dot x^j+N_{ij}\dot x^i\dot x^j.
\]
At the minimum point $(p,v)$, the identities \eqref{eq:xderiv} and \eqref{eq:identderiv}, the vanishing of the connection coefficients, and the
inequality \eqref{eq:minineq} imply the following:
\begin{align*}
\frac{\partial Z}{\partial t}& =a^{kl}\frac{\partial^2S_{11}}{\partial
x^k\partial x^l} +N_{11}\\
&\geq 2a^{kl}\left(2\Gamma_k^p\nabla_lS_{1p}-\Gamma_k^p\Gamma_l^qS_{pq}
\right)+N_{11}.
\end{align*}
The assumption of the theorem implies that the supremum of the right-hand side over all choices of $\Gamma$ is non-negative, so the left-hand side is non-negative.
By the maximum principle (see the argument of Hamilton in
\cite{Ha2}) the inequality $Z\geq 0$ is preserved.
\end{proof}

To illustrate the application of Theorem \ref{thm:newMP} I include the 
following result on preserving convexity for solutions of fully nonlinear
parabolic equations.  A result of this kind was first shown to me by
Gerhard Huisken (in the context of preserving convexity for evolving
hypersurfaces), who proved it by considering the evolution equation for
the inverse of the second derivative matrix.  It can also be proved by
considering the equation satisfied by the second derivatives of
the Legendre transform of the solution.  The belief that the conclusion
should also follow directly from the evolution equation for $D^2u$ led to
Theorem \ref{thm:newMP}.

\begin{theorem}\label{thm:presconv}
Let $\Omega$ be a bounded domain in ${\mathbb R}^n$.  Let 
$u: \Omega\times[0,T]\to{\mathbb R}$ be a solution of a fully nonlinear
equation of the form
\[
\frac{\partial u}{\partial t} = F(D^2u)
\]
where $F$ is a $C^2$ function defined on the cone $\Gamma_+$ of positive 
definite symmetric matrices, which is monotone increasing (that is,
$F(A+B)\geq F(A)$ whenever $B$ is a positive definite matrix), and such
that the function
\[
F^*(A) = -F(A^{-1})
\]
is concave on $\Gamma_+$.  If $D^2u\geq \varepsilon I$ (for some 
$\varepsilon>0$) everywhere on $\Omega$ for $t=0$, and on
$\partial\Omega$ for $0\leq t\leq T$, then $D^2u\geq \varepsilon I$
everywhere on $\Omega$ for $0\leq t\leq T$.
\end{theorem}

\begin{proof}
$D^2u$ evolves as follows (denoting derivatives by subscripts):
\[
\frac{\partial}{\partial t}u_{ij}
= \dot F^{kl}D_kD_lu_{ij}
+ \ddot F^{kl,pq} u_{ikl}u_{jpq}.
\]
The result is obvious for convex $F$, but not for the weaker condition
of the Theorem.  Theorem
\ref{thm:newMP} with 
$S_{ij}=u_{ij}-\varepsilon\delta_{ij}$ and 
$\Lambda_k^p=v^jD_jS_{kq}{\mathfrak r}^{qp}$ (where ${\mathfrak r}$ is
the inverse matrix of $D^2u$) gives
\begin{align*}
\sup_{\Gamma}\left(4a^{kl}\Gamma_k^pD_lS_{ip}v^i
-2a^{kl}\Gamma_k^p\Gamma_l^q
S_{pq}\right)&=
4a^{kl}\Gamma_k^pD_lS_{ip}v^i
-2a^{kl}\Gamma_k^p\Gamma_l^q
S_{pq}\\
&= 2v^iv^j\dot F^{kl}u_{ikq}{\mathfrak r}^{qp}u_{jpl}\\
&\quad\null+2\varepsilon v^iv^j\dot F^{kl}u_{ika}{\mathfrak
r}^{ap}\delta_{pq} {\mathfrak r}^{qb}u_{jlb}\\
&\geq 2v^iv^j\dot F^{kl}u_{ikq}{\mathfrak r}^{qp}u_{jlp}.
\end{align*}
The result then follows from Theorem \ref{thm:newMP} and the following
Lemma.
\end{proof}

\begin{lemma}\label{lem:dualconc}
If $F^*$ is concave, then 
\begin{equation}\label{eq:dualconc}
\left(\ddot F^{kl,pq}+2\dot F^{kp}{\mathfrak r}^{lq}
\right)X_{kl}X_{pq}\geq 0
\end{equation}
for every symmetric matrix $X$.
\end{lemma}

\begin{proof}
Take the identity
$F(A) = -F^*(B)$ with $B=A^{-1}$, and differentiate:
\[
\dot F^{kl}(A) = \frac{\partial F^*}{\partial B^{ab}}(B)B^{ak}B^{bl}.
\]
Further differentiation gives
\begin{align*}
\ddot F^{kl,pq} &= -\frac{\partial^2 F^*}{\partial B^{ab}\partial
B^{cd}} B^{ak}B^{bl}B^{cp}B^{dq}-\frac{\partial F^*}{\partial
B^{ab}}
\left(B^{ak}B^{bp}B^{ql}+B^{ap}B^{kq}B^{bl}\right)\\
&\geq -\frac{\partial F^*}{\partial
B^{ab}}
\left(B^{ak}B^{bp}B^{ql}+B^{aq}B^{kp}B^{bl}\right)\\
&=-{\dot F^{kp}}(A^{-1})^{ql}-{\dot F^{ql}}(A^{-1})^{kp}.
\end{align*}
This proves the Lemma.
\end{proof}

Theorem \ref{thm:presconv} extends (with trivial modifications) to 
equations of the form $\frac{\partial u}{\partial t} = F(D^2u,Du,u)$
if the same concavity condition holds in the first argument,
and $F$ is also convex in the last entry.

\section{The pinching estimate}

The main result of this paper gives conditions under which an equation
will preserve uniform positivity of the second derivatives, in the
sense that $D^2u\geq
\varepsilon
\Delta u\, \text{\rm I}$ for some $\varepsilon\in(0,1/n)$.  This is a
non-trivial extension of Theorem \ref{thm:presconv}, and requires
considerably more work to prove.  The result is stated to
allow easy application in a variety of different situations.  In
the model case Equation \eqref{eq:sysd2u} implies that the tensor
$S_{ij}= D_iD_ju-\varepsilon\Delta u\,\delta_{ij}$ evolves according to
\begin{equation}\label{eq:evoldiff}
\frac{\partial S_{ij}}{\partial t} = \dot F^{kl}D_kD_lS_{ij}
+\ddot F^{kl,pq}u_{ikl}u_{jpq}-\varepsilon
\delta_{ij}\delta^{ab}\ddot F^{kl,pq}u_{akl}u_{bpq}.
\end{equation}
The import of the inequality stated in the Theorem below should be understood
in view of Equation \eqref{eq:evoldiff} and Theorem \ref{thm:newMP}.  In this 
context the tensor $T_{ijk}$ which appears there represents $u_{ijk}$.

\begin{theorem}\label{thm:main}
Let $f$ be a smooth symmetric, monotone, concave and inverse-concave function on 
$\Gamma_+$.  Let $F(A)=f(\lambda(A))$, where
$\lambda$ is the map which takes a symmetric matrix to its eigenvalues.
Let $A$ be a symmetric positive definite matrix and $v$ an
eigenvector of $A$ corresponding to the smallest eigenvalue of $A$, and
let $\varepsilon = \frac{A_{ij}v^iv^j}{\text{\rm Tr A}|v|^2}\in (0,1/n)$.
If $T$ is a totally symmetric $3$-tensor
with $T_{ijk}v^iv^j=\varepsilon \delta^{ij}T_{ijk}$ for $k=1,\dots,n$, then
\begin{align*}
&v^iv^j\ddot F^{kl,pq}(A)T_{ikl}T_{jpq}-\varepsilon|v|^2\delta^{ab}
\ddot F^{kl,pq}(A)T_{akl}T_{bpq}\\
&+2\sup_{\Gamma}\dot
F^{kl}(A)\left(2\Gamma_k^p\left(T_{lpi}v^i-\varepsilon \delta^{ab}T_{lab}v_p\right)-
\Gamma_k^p\Gamma_l^q(A_{pq}-\varepsilon\text{\rm Tr}A\delta_{pq})
\right)\geq 0.
\end{align*}
\end{theorem}

\section{Differentiating eigenvalues and eigenvectors}

This section prepares for the proof of Theorem \ref{thm:main} by establishing
results about derivatives of eigenvalues and eigenvectors of symmetric matrices, 
and of functions of symmetric matrices defined in terms of their eigenvalues.

\begin{theorem}\label{thm:difffn}
Let $f$ be a $C^2$ symmetric function defined on a symmetric region $\Omega$ in 
${\mathbb R}^n$.  Let $\tilde\Omega=\{A\in\text{\rm Sym}(n):\ \lambda(A)\in\Omega\}$,
and define $F: \tilde\Omega\to\Rls$ by $F(A)=f(\lambda(A))$.  Then at any diagonal
$A\in\tilde\Omega$ with distinct eigenvalues, the second derivative of $F$ in
direction $B\in\text{\rm Sym}(n)$ is given by
\[
\ddot F(B,B) = \sum_{k,l}\ddot f^{kl}B_{kk}B_{ll} + 2\sum_{k<l}
\frac{\dot f^k-\dot f^l}{\lambda_k-\lambda_l}B_{kl}^2.
\]
\end{theorem}

This result appeared without detailed proof in \cite[Equation 2.23]{A1}.  A proof 
appeared later in \cite{Ge2}, involving somewhat laborious checking of several
cases.  For this reason I include here an argument which may be more illuminating.

\begin{proof}
Let $Z: \text{\rm Sym}(n)\times{\mathbb R}^n\times O(n)\to\text{\rm Sym}(n)$ be the 
smooth map given by
\[
Z(A,\lambda,M) = M^tAM-\text{\rm diag}(\lambda).
\]
This vanishes if and only if the columns of $M$ are eigenvectors of $A$, with 
eigenvalues $\lambda$.  The derivative of $Z$ at a point $(A,\lambda,M)$
with $Z=0$ in a direction $(A',\lambda',M')$ is as follows (noting that 
$M'=M\Lambda$ with $\Lambda_{ij}+\Lambda_{ji}=0$):
\begin{align}\label{eq:DZ}
DZ(A', \lambda',M')_{ij}&= M_{ik}A'_{kl}M_{lj}+
 \Lambda_{ki}M_{lk}A_{lp}M_{pj}+M_{ki}A_{kl}M_{lp}\Lambda_{pj}
-\lambda'_i\delta_{ij}\notag\\
&=M_{ik}A'_{kl}M_{lj} + \left(\lambda_j\Lambda_{ji}+\lambda_i\Lambda_{ij}\right)
-\lambda'_i\delta_{ij}\notag\\
&=M_{ik}A'_{kl}M_{lj} +(\lambda_i-\lambda_j)\Lambda_{ij}-\lambda'_i\delta_{ij}.
\end{align}
The restriction of this to the last two components has no kernel:  If it vanishes, 
the diagonal parts imply $\lambda'=0$, and the off-diagonal parts imply $\Lambda=0$
since $\lambda_i-\lambda_j\neq 0$.  Therefore this is an isomorphism, and the
implicit function theorem gives that the zero set of $Z$ is locally of the form
$\{\lambda=\lambda(A),\ M=M(A)\}$ where $\lambda$ and $M$ are analytic functions of
$A$.  If $A$ is diagonal (so $M=I$) the first derivatives of $\lambda$ and $M$ can be
read off:
\[
\lambda'_i = A'_{ii}
\qquad\text{\rm and}\qquad
M'_{ij} = \Lambda_{ij} = -\frac{A'_{ij}}{\lambda_i-\lambda_j}.
\]
Equation \eqref{eq:DZ} holds everywhere on $\{Z=0\}$, so differentiating this along 
$(A(t),\lambda(A(t)),M(A(t))$ with $A''=0$ and $M(0)=I$ gives
\[
0 = M'_{ki}A'_{kj}+A'_{ik}M'_{kj}+(\lambda'_i-\lambda'_j)
\Lambda_{ij}+(\lambda_i-\lambda_j)\Lambda'_{ij}-\lambda''_i\delta_{ij}.
\]
The second derivative of $\lambda_i$ can be read off from the $(i,i)$ component:
\[
\lambda''_i = \sum_k\left(\Lambda_{ki}A'_{ki}+A'_{ik}\Lambda_{ki}\right)
=-2\sum_{k\neq i}\frac{\left(A'_{ik}\right)^2}{\lambda_k-\lambda_i}.
\]
The first and second derivatives of $F$ at a diagonal matrix $A$ with distinct 
eigenvalues can now be computed directly:
\begin{equation}\label{eq:dotf}
F'= \sum_k\dot f^k\lambda'_k
=\sum_k\dot f^kA'_{kk}
=\sum_{k,l}\dot f^k\delta_{kl}A'_{kl},
\end{equation}
so $\dot F^{kl}=\dot f^k\delta_{kl}$.  Also,
\begin{align*}
F'' &= 
\frac{d}{dt}\left(\sum_k\dot f^k\lambda'_k\right)\\
&=\sum_{k,l}\ddot f^{kl}\lambda'_k\lambda'_l
+\sum_k\dot f^k\lambda''_k\\
&=\sum_{k,l}\ddot f^{kl}A'_{kk}A'_{ll}
-2\sum_{k\neq l}\frac{\dot f^k}{\lambda_l-\lambda_k}\left(A'_{kl}\right)^2\\
&=\sum_{k,l}\ddot f^{kl}A'_{kk}A'_{ll}
+2\sum_{k<l}\frac{\dot f^k-\dot f^l}{\lambda_k-\lambda_l}\left(A'_{kl}\right)^2.
\end{align*}
This completes the proof.
\end{proof}

\begin{corollary}\label{cor:conc}
$F$ is concave at $A$ if and only if $f$ is concave at $\lambda(A)$ and 
\[
\frac{\dot f^k-\dot f^l}{\lambda_k-\lambda_l}\leq 0\qquad\text{\rm for all } k\neq l.
\]
\end{corollary}

\begin{corollary}\label{cor:concequiv}
For $\Omega$ convex, $F$ is concave on $\tilde\Omega$ if $f$ is concave on $\Omega$.
\end{corollary}

\begin{proof}
See \cite[Lemma 2.2]{A1}.
\end{proof}

\begin{corollary}\label{cor:invconc}
$F^*$ is concave at $A$ if and only if 
\[
\left(\ddot f^{kl}+2\frac{\dot f^k}{\lambda_k}\delta_{kl}\right)
\geq 0\qquad\text{\rm and}\qquad 
\frac{\dot f^k-\dot f^l}{\lambda_k-\lambda_l}+\frac{\dot f^k}{\lambda_l}
+\frac{\dot f^l}{\lambda_k}\geq 0,\quad k\neq l.
\]
\end{corollary}

\begin{proof}
This follows immediately from the inequality \eqref{eq:dualconc} and Theorem 
\ref{thm:difffn}.
\end{proof}

\begin{corollary}\label{cor:dualconcequiv}
If $\Omega^{-1} = \{(x_1^{-1},\dots,x_n^{-1}):\ (x_1,\dots,x_n)\in\Omega\}$ is 
convex,
then $F^*$ is concave if and only if $f^*$ is concave.
\end{corollary}

In particular, corollaries \ref{cor:concequiv} and \ref{cor:dualconcequiv} apply
for functions defined on $\Gamma_+$.

\section{Proof of the estimate}

This section contains the proof of the main result, Theorem \ref{thm:main}.  
Note that the Theorem does not refer at all to a partial differential equation 
or its solution, but only to a pointwise inequality for the first and second 
derivatives of a function defined on the positive cone.  If $F$ is $C^2$, then $\ddot
F$ and $\dot F$ are continuous, so for fixed $\Gamma$, $v$ and $T$ the quantity we
wish to estimate, 
\begin{align*}
&v^iv^j\ddot F^{kl,pq}(A)T_{ikl}T_{jpq}-\varepsilon|v|^2\delta^{ab}
\ddot F^{kl,pq}(A)T_{akl}T_{bpq}\\
&+2\dot
F^{kl}(A)\left(2\Gamma_k^p\left(T_{lpi}v^i-\varepsilon T_{laa}v_p\right)-
\Gamma_k^p\Gamma_l^q(A_{pq}-\varepsilon\text{\rm Tr}A\delta_{pq})
\right)
\end{align*}
is continuous in $A$.  It follows that the supremum over $\Gamma$ is 
semi-continuous in $A$.  We will take advantage of this by working only with 
symmetric
matrices $A$ for which all of the eigenvalues are distinct.  This is possible 
since
for any positive definite $A\in\text{\t\rm Sym}(n)$ with 
$A_{ij}\geq\varepsilon\text{\rm Tr}A\delta_{ij}$ and $A_{ij}v^iv^j =
\varepsilon\text{\rm Tr}A|v|^2$ for some $v\neq 0$, there is a sequence
$\{A^{(k)}\}_{k\geq 0}$ approaching $A$,
satisfying $A^{(k)}_{ij}\geq\varepsilon\text{\rm Tr}A^{(k)}\delta_{ij}$ and
$A^{(k)}_{ij}v^iv^j = \varepsilon\text{\rm Tr}A^{(k)}|v|^2$, and with each
$A^{(k)}$ having distinct eigenvalues.  Hence it suffices to establish the result in
the case where all of the eigenvalues are distinct.

In this case there is an orthonormal basis
$e_1,\dots,e_n$ consisting of eigenfunctions of $A$, with eigenvalues in increasing 
order.  In this basis, $v=e_1$ and 
$A=\text{\rm diag}(\lambda_1,\dots,\lambda_n)$, and
$\lambda_1=\varepsilon(\lambda_1+\dots+\lambda_n)$.   Also 
$\dot F=\text{\rm diag}(\dot f^1,\dots,\dot f^n)$ by Equation \eqref{eq:dotf}.

The problem is simplified by the observation that the supremum over 
$\Gamma$ can be computed exactly in this case:  We can write
\begin{align*}
&2\dot F^{kl}\left(2\Gamma_k^p\left(T_{lpi}v^i-\varepsilon T_{laa}v_p\right)-
\Gamma_k^p\Gamma_l^q(A_{pq}-\varepsilon\text{\rm Tr}A\delta_{pq})\right)\\
&= 2\sum_{k=1}^n\sum_{p=2}^n\dot f^k\left(\Gamma_k^pT_{kp1}-(\Gamma_k^p)^2
(\lambda_p-\lambda_1)\right)\\
&=2\sum_{k\geq 1,p\geq 2}\left(\frac{\dot f^k}{\lambda_p-\lambda_1}T_{kp1}^2
-\dot f^k(\lambda_p-\lambda_1)
\left(\Gamma_k^p-\frac{T_{kp1}}{\lambda_p-\lambda_1}\right)^2\right)
\end{align*}
It follows that the supremum is attained by the choice 
$\Gamma_k^p = \frac{T_{kp1}}{\lambda_p-\lambda_1}$.

The required inequality becomes the following:
\begin{align*}
Q &= \sum_{k,l}\ddot f^{kl}T_{1kk}T_{1ll}
-\varepsilon\sum_{j,k,l}\ddot f^{kl}T_{jkk}T_{jll}
+2\sum_k\sum_{l>1}\frac{\dot f^k}{\lambda_l-\lambda_1}T_{1kl}^2\\
&\quad\null + 2\sum_{k< l}\frac{\dot f^k-\dot f^l}{\lambda_k-\lambda_l}T_{1kl}^2
-2\varepsilon\sum_j\sum_{k<l}\frac{\dot f^k-\dot f^l}{\lambda_k-\lambda_l}
T_{jkl}^2.\\
&\geq 0.
\end{align*}

We use the identities
$$
T_{k11} = \frac{\varepsilon}{1-\varepsilon}\sum_{j>1}T_{kjj}
$$
for $k=1,\dots,n$, to eliminate terms involving $T_{k11}$.  This together with 
the total symmetry of $T$ implies that, as a bilinear form on the space of all
possible $T$, $Q$ has a block-diagonal form, as follows:
\[
Q = \sum_{k=1}^nQ_k+\sum_{1\leq j<k<l}Q_{jkl}
\]
where $Q_k$ involves only $T_{kii}$ for $i\geq 2$, and $Q_{jkl}$ involves only
$T_{jkl}$.  Precisely, these are as follows:
\begin{align*}
Q_1& = (1-\varepsilon)\sum_{k,l>1}\left(\dot f^{kl} +
\frac{\varepsilon}{1-\varepsilon}(\ddot f^{k1}+\ddot f^{1l})+\left(
\frac{\varepsilon}{1-\varepsilon}\right)^2\ddot f^{11}\right)T_{1kk}T_{1ll}\\
&\quad\null + 2\sum_{k>1}\frac{(1-\varepsilon)\dot f^k+\varepsilon \dot
f^1}{\lambda_k-\lambda_1}T_{1kk}^2;\\
Q_k&=-\varepsilon\sum_{i,j>1}\left(\ddot
f^{ij}+\frac{\varepsilon}{1-\varepsilon}(\ddot f^{i1}+\ddot f^{1j})
+\left(\frac{\varepsilon}{1-\varepsilon}\right)^2\ddot
f^{11}\right)T_{kii}T_{kjj}\\
&\quad\null +2\frac{(1-\varepsilon)\dot f^k+\varepsilon \dot f^1}
{\lambda_k-\lambda_1}\left(\frac{\varepsilon}{1-\varepsilon}\sum_{i>1}T_{kii}
\right)^2
-2\varepsilon\sum_{j\neq k,1}\frac{\dot f^k-\dot f^j}{\lambda_k-\lambda_j}
T_{kjj}^2;\\
Q_{1kl}&=2\left((1-\varepsilon)\frac{\dot f^k-\dot
f^l}{\lambda_k-\lambda_l}+\frac{\dot f^l}{\lambda_k-\lambda_1}+\frac{\dot
f^k}{\lambda_l-\lambda_1}-\varepsilon\frac{\dot
f^k-\dot
f^1}{\lambda_k-\lambda_1}-\varepsilon\frac{\dot f^l-\dot
f^1}{\lambda_l-\lambda_1}\right)T_{1kl}^2;\\ 
Q_{jkl}&=-2\varepsilon
\left(\frac{\dot f^k-\dot f^l}{\lambda_k-\lambda_l}
+\frac{\dot f^k-\dot f^j}{\lambda_k-\lambda_j}
+\frac{\dot f^l-\dot f^j}{\lambda_l-\lambda_j}\right)T_{jkl}^2.
\end{align*}

We require each of these to be non-negative.

For $k>l>j>1$, $Q_{jkl}\geq 0$ by concavity (see Corollary \ref{cor:conc}). 

The same is true for $Q_k$:  The matrix in the first bracket is
$\ddot f(\xi,\xi)$, where
\begin{equation*}
\xi = \sum_{k>2}T_{1kk}\left(e_k+\frac{\varepsilon}{1-\varepsilon}e_1\right).
\end{equation*}
The concavity of $f$ therefore implies that this term is non-positive.
The last term is also of the right sign
by Corollary \ref{cor:conc}, and the remaining term is manifestly non-negative.

The non-negativity of $Q_{1kl}$ follows from the
concavity of both $f$ and $f^*$:
\begin{align*}
Q_{1kl}&=2(1-\varepsilon)\left(\frac{\dot f^k-\dot f^l}{\lambda_k-\lambda_l}
+\frac{\dot f^k}{\lambda_l-\lambda_1}+\frac{\dot f^l}{\lambda_k-\lambda_1}
\right)\\
&\quad\null + 2\varepsilon\left(\frac{\dot f^k}{\lambda_l-\lambda_1}
+\frac{\dot f^l}{\lambda_k-\lambda_1}-
\frac{\dot f^k-\dot f^1}{\lambda_k-\lambda_1}-
\frac{\dot f^l-\dot f^1}{\lambda_l-\lambda_1}\right).
\end{align*}
The first bracket is non-negative by the second inequality of Corollary 
\ref{cor:invconc}.  The first two terms in the second bracket are manifestly
non-negative, while the other two are non-negative by Corollary \ref{cor:conc}.

Finally, non-negativity of $Q_1$ follows from concavity of $f^*$ in an
indirect way:  Consider the function $\phi$ of $(n-1)$ variables $x_2,\dots,x_n$ 
given by
$\phi(x_2,\dots,x_n)=f\left(\frac{\varepsilon}{1-\varepsilon}
(x_2+\dots+x_n),x_2,\dots,x_n
\right)$.
Then
\begin{align*}
Q_1 &= (1-\varepsilon)\sum_{k,l>1}\left(\ddot\phi^{kl}+2\frac{\dot\phi^k}
{\lambda_k-\lambda_1}\delta_{kl}
\right)T_{1kk}T_{1ll}\\
&\geq 
(1-\varepsilon)\sum_{k,l>1}\left(\ddot\phi^{kl}+2\frac{\dot\phi^k}{\lambda_k}
\delta_{kl}
\right)T_{1kk}T_{1ll}.
\end{align*}
The first inequality in Corollary \ref{cor:invconc} then implies
$Q_1\geq 0$ provided $\phi^*$ is concave.  To establish concavity of $\phi^*$, note
\begin{align*}
\phi^*(x_2,\dots,x_n)&=
-\phi\left(\frac{\varepsilon}{1-\varepsilon}(x_2^{-1}+\dots+x_n^{-1}),x_2^{-1},
\dots,x_n^{-1}\right)\\
&=f^*(\psi(x_2,\dots,x_n),x_2,\dots,x_n),
\end{align*}
where $\psi(x_2,\dots,x_n) = \frac{1-\varepsilon}{\varepsilon}
\left(x_2^{-1}+\dots+x_n^{-1}\right)^{-1}$ is a multiple of the harmonic 
mean of
$x_2,\dots,x_n$, hence a concave function of $x_2,\dots,x_n$ (Corollary \ref{cor:power}).  Also, $f^*$ is concave, and increasing in each argument.  Therefore
\begin{align*}
\alpha\phi^*(x)+(1-\alpha)\phi^*(y)&=\alpha f^*(\psi(x),x)+(1-\alpha)f^*(\psi(y),y)\\
&\leq f^*(\alpha\psi(x)+(1-\alpha)\psi(y),\alpha x+(1-\alpha)y)\\
&\leq f^*(\psi(\alpha x+(1-\alpha)y),\alpha x+(1-\alpha)y)\\
&=\phi^*(\alpha x+(1-\alpha)y)
\end{align*}
for any $\alpha\in(0,1)$ and $x,y$ in the positive cone of 
${\mathbb R}^{n-1}$.
Here the first inequality follows from the concavity of $f^*$, and the second 
follows since $\psi$ is concave (so $\alpha\psi(x)+(1-\alpha)\psi(y)\leq \psi(\alpha
x+(1-\alpha)y)$) and $f^*$ is increasing in the first argument. Therefore $\phi^*$ is
concave, and the proof is complete.

\section{Application to evolving hypersurfaces}

In this section Theorem \ref{thm:main} is applied to prove Theorem
\ref{thm:contraction} on evolving hypersurfaces.  As mentioned before,
the only case not proved elsewhere is case 4, where both the speed
$f$ and its dual $f^*$ are concave functions on $\Gamma_+$.

The only new ingredient in the proof is the
application of Theorem \ref{thm:main} to prove that the smallest
eigenvalue of $\frac{h_{ij}}{H}$ over $M_t$ is non-decreasing in $t$,
where $h_{ij}$ is the second fundamental form and $H$ is the mean
curvature.

First note that in any local coordinates for $M$ the rate of change
of the metric tensor $g_{ij}=g(\partial_i,\partial_j)$ under Equation
\eqref{eq:flow} is given by
\[
\frac{\partial}{\partial t}g_{ij}=-2Fh_{ij}.
\]
The evolution of $h_{ij}$ is as follows (see
\cite[Lemma 3.13]{A1}):
\[
\frac{\partial}{\partial t}h_{ij}=\dot F^{kl}\nabla_k\nabla_lh_{ij}
+\ddot F^{kl,pq}\nabla_ih_{kl}\nabla_jh_{pq}+h_{ij}\dot
F^{kl}h_{kp}g^{pq}h_{ql}-2Fh_{ik}g^{kl}h_{lj}.
\]
Suppose $S_{ij}=h_{ij}-\varepsilon H g_{ij}$ is non-negative.  Then note
\begin{align}\label{eq:evolS}
\frac{\partial}{\partial t}S_{ij} &= 
\dot F^{kl}\nabla_k\nabla_lS_{ij}
+\ddot F^{kl,pq}T_{ikl}T_{jmn}-\varepsilon g_{ij}
g^{rs}\ddot F^{kl,pq}T_{rkl}T_{spq}\\
&\quad\null+S_{ij}\dot
F^{kl}h_{kp}g^{pq}h_{ql}-2Fh_{ik}g^{kl}S_{lj}\notag
\end{align}
where $T_{ijk} = \nabla_ih_{jk}$ is totally symmetric by the
Codazzi identity, and satisfies
$0=\nabla_kS_{11}=T_{k11}-\varepsilon\sum_jT_{kjj}$ at a point where $S_{ij}$ has
$e_1$ as a null eigenvector.  Then
Theorem
\ref{thm:main} gives 
\begin{align*}
0&\leq\ddot F^{kl,pq}T_{1kl}T_{1mn}-\varepsilon g^{rs}
\ddot F^{kl,pq}T_{rkl}T_{spq}\\
&\quad\null+2\sup_{\Gamma}\dot
F^{kl}\left(2\Gamma_k^p\left(T_{lp1}-\varepsilon T_{ljj}\delta_{p1}\right)-
\Gamma_k^p\Gamma_l^qS_{pq}
\right)
\end{align*}
The terms in the second line of Equation \eqref{eq:evolS} vanish at a null eigenvector, so Theorem \ref{thm:newMP} implies that
$S_{ij}$ remains non-negative.

The remainder of the proof of Theorem
\ref{thm:contraction} is the same as in \cite{A1}.



\begin{thebibliography}{99}

\small

\bibitem{A1} B.~Andrews, {\it Contraction of convex hypersurfaces in Euclidean 
space}, Calc. Var. P.D.E. {\bf 2} (1994), 151--171.

\bibitem{A2} \bysame, {\it Fully nonlinear parabolic equations in two space variables}, preprint, 8 pages, arXiv: math.AP/0402235

\bibitem{A3} \bysame, {\it Moving surfaces by non-concave curvature functions}, 
preprint, 8 pages, arXiv: math.DG/0402273

\bibitem{Lb} G.~Lieberman, ``Second order parabolic differential equations'', World
Scientific, 1996.

\bibitem{Ch1} B.~Chow, {\it Deforming convex hypersurfaces by the $n$th root 
of the Gaussian curvature}, J. Differential Geometry {\bf 23} (1985), 117--138.

\bibitem{Ch2} \bysame, {\it Deforming hypersurfaces by the square root 
of the scalar curvature}, Invent. Math. {\bf 87} (1987), 63--82.

\bibitem{Ge2} C. Gerhardt, {\it Closed Weingarten hypersurfaces in Riemannian
manifolds}, J. Differential Geom. {\b 43} (1996), 612--641.

\bibitem{Ha1} R.~Hamilton, {\it Three-manifolds with positive Ricci curvature}, J. Differential Geometry {\bf 17} (1982), 255--306.

\bibitem{Ha2} \bysame, {\it
Four-manifolds with positive curvature operator}, 
J. Differential Geometry {\bf 24} (1986), 153--179.

\bibitem{Hu1} G.~Huisken, {\it Flow by mean curvature of convex hypersurfaces into 
spheres}, J. Differential Geometry {\bf 20} (1984), 237--268.


\end{thebibliography}
\end{document}